\documentclass[11pt,a4paper,twoside]{article}
\usepackage[centertags]{amsmath}
\usepackage{amsfonts}
\usepackage{amssymb}
\usepackage{amsthm}
\usepackage{latexsym}
\usepackage{newlfont}
\usepackage{verbatim}
\usepackage[english]{babel}
\usepackage[arrow, matrix]{xy}
\theoremstyle{plain}
\newtheorem{thm}{Theorem}
\newtheorem{cor}[thm]{Corollary}
\newtheorem{lem}[thm]{Lemma}
\newtheorem{prop}[thm]{Proposition}

\theoremstyle{definition}
\newtheorem{defn}{Definition}

\newcommand{\Log}{\mathsf{Log}}
\newcommand{\Top}{\mathsf{Top}}
\newcommand{\eq}{\leftrightsquigarrow}
\newcommand{\N}{\mathbb{N}}

\newcommand{\Ext}[2]{\ensuremath{[\![#2]\!]^{#1}}}

\newcommand{\Abelard}{$\forall$belard }
\newcommand{\Eloise}{$\exists$loise }
\newcommand{\C}{\mathcal{C}}
\newcommand{\Q}{\mathbb{Q}}
\newcommand{\noip}{\Log(T_1^{\textrm{no i.p.}})}
\newcommand{\T}{\mathcal{T}}

\begin{document}


\bibliographystyle{amsplain} 

\author{Dmitry Sustretov}
\title{Topological Semantics and Decidability}

\maketitle

\begin{abstract}
  It is well-known that the basic modal logic of all topological
  spaces is $S4$. However, the structure of basic hybrid logics of
  classes of spaces satisfying various separation axioms was until
  present unclear.  We give a direct proof of that modal logics of
  $T_0$, $T_1$ and $T_2$ topological spaces coincide and are $S4$. We
  also examine basic hybrid logics of these classes and prove their
  decidability; as part of this, we find out that the hybrid logics of
  $T_1$ and $T_2$ spaces coincide. Finally, we prove that logics of
  $T_0$ and $T_1$ spaces are PSPACE-complete.
\end{abstract}

\section{Basic definitions}

In this paper we are going to study modal logics that arise as sets of
all formulas valid on certain classes of topological spaces. Thus the
first definition in this paper is bound to be about how the modal
formulas are interpreted on topological spaces (topological semantics
was first introduced by Tarski \cite{tarski}).

\begin{defn}[Topological semantics] 
  A topological space is a pair $(T, \tau)$ where $\tau \subseteq
  \mathcal{P}(T)$ such that $\emptyset, T \in \tau$ and $\tau$ is
  closed under finite intersections and arbitrary unions. Elements of
  $\tau$ are called {\em opens}, an open containing a point $x$ is
  called a {\em neighborhood} of the point $x$.

  A {\em topological model} $\mathfrak{M}$ is a tuple $(T,\tau,V)$ where
  $(T,\tau)$ is a topological space and the valuation
  $V:\textsc{Prop}\to\mathcal{P}(T)$ sends propositional letters to
  subsets of $T$.

  Truth of a formula $\phi$ (of the basic modal language) at a point
  $w$ in a topological model $\mathfrak{M}$ (denoted by
  $\mathfrak{M},w \models \phi$) is defined inductively:
\[\begin{array}{lll}
\mathfrak{M},w \models p & \textrm{ iff } & x \in V(p)\\
\mathfrak{M},w \models \phi \land \psi & \textrm{ iff } &
\mathfrak{M},w
\models \phi \textrm{ and } \mathfrak{M},w \models \psi\\
\mathfrak{M},w \models \neg \phi & \textrm{ iff } & \mathfrak{M},w
\nvDash
\phi\\
\mathfrak{M},w \models \Box \phi & \textrm{ iff } & \exists O \in \tau
\textrm{ such that } w\in O \textrm{ and } \forall v\in
O.(\mathfrak{M},v \models \phi)
\end{array}\]
The basic modal language can be extended with nominals and $@$
operator (in this case we call it $H(@)$) and universal modality $A$
(we denote the dual modality $E$ and call the language
$H(E)$). Nominals are a special kind of propositional letters: it is
required that their valuation is a singleton set. The semantics of $@$
and $E$ is given below:

\[\begin{array}{lll}
\mathfrak{M},w \models @_i \varphi & \textrm{ iff } & \exists v\,
\mathfrak{M},v \models i \textrm{ and } \mathfrak{M},v \models
\varphi\\
& & \textrm{(where $i$ is a nominal)}\\
\mathfrak{M},w \models E\varphi & \textrm{ iff } &  \exists v\,
\mathfrak{M},v \models \varphi\\
\end{array}
\]
\end{defn}

Relational and topological semantics are not completely unrelated; it
is possible to transform certain topological spaces into frames and
vice versa in a satisfiability-preserving fashion.

\begin{prop} \label{prop:transform} 
A topological space is called {\em Alexandroff} if every point of
that space has a minimal neighborhood.

For any Alexandroff space $(T, \tau)$ there exists a binary relation
$R$ such that for any valuation $V$ and for any formula $\varphi \in
H(E)$, $(T, R, V), w \models \varphi$ iff $(T, \tau, V), w \models
\varphi$.

For any transitive reflexive frame $(W, R)$ there exists a topology
$\tau$ on $W$ such that for any valuation $V$ and for any formula
$\varphi \in H(E)$, $(W, R, V), w \models \varphi$ iff $(W, \tau, V),
w \models \varphi$.
\end{prop}

\begin{proof}
  See \cite{johan}, section 2.4.
\end{proof}

\begin{defn}[Topobisimulation]
  Let $(T, \tau, V)$ and $(S, \sigma, W)$ be two topological models
  and consider a relation $R \subseteq T \times S$. Denote
  $$
  \begin{array}{l}
    R(X) = \{ y \mid \exists x \in X, (x,y) \in R\}\\
    R^{-1}(Y) = \{ x \mid \exists y \in Y, (x,y) \in R\}
  \end{array}
  $$
  for any subset $X \subseteq T$, $Y \subseteq S$.

  The relation $R$ is called a {\em topobisimulation} if
  
  \begin{itemize}
  \item[\bf Prop] if $Rxy$ then for all $p \in \textsc{Prop}$,
    $(T,\tau),V,x\models p$ iff $(S,\sigma),W,y\models p$  
  \item[\bf Zig] for any $O \in \tau$, $R(O) \in \sigma$
  \item[\bf Zag] for any $U \in \sigma$, $R^{-1}(U) \in \tau$
  \end{itemize}  

  A bisimulation is called {\em total} iff for any $x \in T$ there is
  $y \in S$ such that $Rxy$ and for any $y \in S$ there is $x \in T$
  such that $Rxy$.

  A bisimulation is called an {\em hybrid} if additionally for any
  nominal $i$ if $x \in V(i)$ and $y \in W(i)$ then $Rxy$.

  A map is called {\em interior} if it is open and
  continuous. Clearly, the graph of an interior map is satisfies {\bf
  Zig} and {\bf Zag} conditions.
\end{defn}

In topological semantics just like in the relational semantics, two
points connected by a topobisimulation satisfy the same formulas (if
the topobisimulation is total, this is true for the formulas with
universal modality). See \cite{johan} for the proofs.

It is well-known that the (basic modal) logic of all topological
spaces is $S4$. In what follows, we are going to deal with three
classes of topological spaces, defined by the so-called separation
axioms.

\begin{defn}[Separation axioms]
  \noindent \begin{enumerate} 
  \item[$T_0$] for any two distinct points $x, y$ there is either an
    open neighborhood of $x$ that does not contain $y$, or an open
    neighborhood of $y$ that does not contain $x$.
  \item[$T_1$] any singleton set is closed.
  \item[$T_2$] any two distinct points $x,y$ can be separated by two open
    neighborhoods, i.e.\ there exist $O_x \ni x, O_y \ni y$ such that
    $O_x \cap O_y = \emptyset$.
  \end{enumerate}
  
\end{defn}

There are necessary and sufficient conditions (given in
\cite{tengasus}) of whether a class of spaces is definable in $H(@)$
(and $H(E)$). Thus, axioms $T_0$ and $T_1$ are definable in $H(@)$,
the formulas are, respectively, $@_i\neg j \to (@_i \Box \neg j \lor
@_j \Box \neg i)$ and $\Diamond i \to i$. On the other hand,
\cite{tengasus} show that $T_2$ is not definable even in $H(E)$. The
basic modal language is even less expressible: none of the separation
axioms is definable in it. Nonetheless, although we know the
boundaries of expressivity of modal and hybrid languages, we know very
little about the structure of the logics. Are the logics of separation
axioms distinct? Are they decidable? If yes, what is the complexity?
In this paper we will address all those questions.

\section{Modal logic with universal modality} \label{sect:bml}

The McKinsey-Tarski theorem \cite{mckinsey},\cite{tarski} states that
the logic of every metrizable dense-in-itself (without isolated
points) space is S4. Thus, for example, the logic of rationals is
$S4$. It follows that the logics of all separation axioms that $\Q$
satisfies are all $S4$. Shehtman \cite{shehtman} showed the same kind
of result holds for the basic modal language enriched with universal
modality. Denote by $\Log(K)$ a set of formulas in basic modal
language with universal modality which are valid on all topological
spaces in class $K$. In this section we will give a direct proof of
these facts, namely that $\Log(T_0)=\Log(T_1)=\Log(T_2)=\Log(\Top)$
where $\Top$ is the class of all topological spaces. Our technique
will be to build a topobisimulation between a finite topological space
and a space from each respective class.

\begin{defn}[Finite model property]
  A logic $L$ has {\em finite model property} with respect to a class
  of topological models $K$ if $K \models L$ and for any $\varphi
  \notin L$ there exists a finite $\mathfrak{M} \in K$ such that
  $\varphi$ is satisfiable on $\mathfrak M$.
\end{defn}

\begin{prop} \label{prop:fmp}
  The logic $S4$ has a finite model property.
\end{prop}

\begin{proof}
  The proof is Proposition~\ref{prop:transform} combined with the fact
  that $S4$ has a finite model property with respect to transitive
  reflexive frames.
\end{proof}

\begin{thm} \label{thm:t0}
$\Log(T_0)=\Log(\Top)$.
\end{thm}

\begin{proof}
  The inclusion $S4 \subseteq \Log(T_0)$ is obvious, so we only have
  to prove $\Log(T_0) \subseteq S4$. Take an arbitrary topological
  space $(T, \tau)$ and define an equivalence relation: $x \eq y$ iff
  for all $O \in \tau, x \in O \textrm{ iff } y \in O$. The quotient
  set with the maximal topology that makes the natural projection
  continuous (this topological space is known as Kolmogorov quotient
  of $T$) is a $T_0$ space. The graph of the natural projection map is
  a topobisimulation. It follows that every formula, that is not an
  $S4$ validity is not a $T_0$ validity either.
\end{proof}

\begin{thm} \label{thm:t1}
Any finite transitive and reflexive model is an interior image of a
$T_1$ space.
\end{thm}

\begin{proof}
  Let $(T, \tau, V)$ be a finite topological model, let us construct a
  topobisimilar model $(S, \sigma, W)$ where $(S, \sigma)$ is a $T_1$
  topological space.

  We will identify $T$ with the initial segment of natural numbers, so
  $T=\{1,\ldots,n\}$. First, let us introduce some notation:
  $$
  X_k = \{ nx + k \mid 0 \leq x < \infty \}, 1 \leq k \leq n
  $$

  Let $\sigma_0$ be a cofinite topology on $\N$,
  that is 
  $$
  \sigma_0 = \{ O \mid \N - O \textrm{ is finite} \}
  $$

  and for any subset $O \subseteq T$ denote

  $$
  \tilde{O}=\{X_i \mid i \in O \}
  $$

  Then define the topology $\sigma$ on $S=\N$ to be generated by the
  set 
  $$
  \beta=\sigma_0 \cup \{ \tilde{O} \mid O \in \tau \}
  $$

  Define valuation to be
  $$
  W(p)=\widetilde{V(p)} \textrm{ for all } p\in \textsc{Prop}
  $$

  $(S,\sigma)$ is a $T_1$ space, because $\sigma$ contains $\sigma_0$,
  hence all complements of singleton sets are open.

  Let $f: S \to T$ be a function that maps $X_k$ to $k$ for all $k \in
  T$. It is an onto interior map.

  First, let us prove that it is open. The topology generated by
  $\beta$ consists of sets of the form $\tilde{O} \cap U$, where $U$
  is an open set from $\sigma_0$ and $O \in \tau$. We have
  $f(\tilde{O} \cap U) = f(\cup_{k \in O} (X_k \cap U))$ and since
  $X_k$ are dense sets in the topology $\sigma_0$, all intersections
  $X_k \cap U$ are non-empty, hence $f(X_k \cap U)=\{k\}$ and
  $f(\tilde{O} \cap U) = \cup_{k \in O} \{k\} = O$. That is, any open
  set is mapped by $f$ to an open set.

  On the other hand, $f$ is continuous. The $f$-preimage of any open
  set $O \in \tau$ is $\tilde{O}$ by construction of $f$, and
  $\tilde{O}$ is open by definition of the topology $\sigma$.

  A final note: consider the graph $R$ of $f$. The points connected by
  $R$ agree on propositional letters, hence $R$ is a
  topobisimulation. Since $f$ is defined on the whole $S$ and is onto,
  it is a total topobisimulation.
  
\end{proof}

\begin{cor}
 $\Log(T_1)=\Log(\Top)$
\end{cor}
\begin{proof}
  By Proposition~\ref{prop:fmp}, $\Log(\Top)$ has finite (topological)
  model property, i.e.\ every formula that is not valid on the class
  of all topological spaces can be refuted on a finite topological
  model. By the theorem we have just proved there exists a total
  bisimulation between any finite topological model and a model based
  on a topological space with $T_1$ topology. It follows that any
  $\Log(\Top)$ non-theorem can be refuted on a $T_1$ model, hence
  $\Log(T_1) \subseteq \Log(\Top)$, hence $\Log(T_1)=\Log(\Top)$.
\end{proof}

The main idea of the proof is that we could construct a space which is
a disjoint union of dense subsets and then add necessary opens to the
topology to get a space topobisimilar to $T$. We can exploit this idea
in a more general setting, leading us to the following

\begin{thm} \label{thm:t2}
Any finite transitive and reflexive model is an interior image of a
$T_2$ space.
\end{thm}

\begin{proof}
  We will use a construction by L.\ Feng and O.\ Masaveu. In the paper
  \cite{feng} they prove that for any cardinal $\alpha$ there exists a
  $T_2$ space which is a disjoint union of $\alpha$ dense subsets
  (such a space is called $\alpha$-resolvable). We will apply this
  statement for a finite $\alpha=n$, so let $(S,\sigma_0)$ be such a
  space and $S=\bigcup_{k=1}^n X_k$ where $X_k$ are disjoint dense
  subsets from the theorem of Feng and Masaveu.

  Now define $\tilde{O}$, $\beta$, $\sigma$, the valuation $W$ and the
  function $f: S \to T$ the same way as in the proof of the previous
  theorem.

  $(S,\sigma)$ stays a $T_2$ space (because we have only added more
  opens to it). For any $O \in \tau$, $f^{-1}(O)=\tilde{O}$ and is
  open by the definition of topology $\sigma$. 

  In order to prove openness of $f$, we will use the same argument as
  in the Theorem~\ref{thm:t1}. Indeed take any open $\tilde{O} \cap U$
  from $\sigma$, where $O \in \tau$ and $U \in \sigma_0$.  Since $X_k$
  are dense sets in $\sigma_0$, all intersections $X_k \cap U$ are
  non-empty, hence $f(X_k \cap U)=O$, which is open.

  Again, we have constructed the valuation in such a way that the
  points connected by the graph $R$ of $f$ agree on propositional
  letters, hence $R$ is a total topobisimulation.
\end{proof} 

\begin{cor}
$\Log(T_2)=\Log(\Top)$
\end{cor}

In fact, nothing in the proof depends on the $T_2$ condition, except
the existence of $n$-resolvable sets. That leads us to the following
more general result.

\begin{thm}
  Let $K$ be a class of topological spaces that contains an
  $n$-resolvable space for all finite $n$ and for any space $(T,\tau)
  \in K$, it is true that $(T,\tau')\in K$ for all $\tau' \supset
  \tau$ (i.e.\ $K$ is closed under refinement of topology). Then
  $\Log(K)=\Log(\Top)$.
\end{thm}

\begin{proof}
  In this theorem we extract the properties of the class $T_2$ used in
  the proof of the Theorem~\ref{thm:t2}. Indeed, we need an
  $n$-resolvable space to start our construction, then we add new
  opens to this space in order to obtain a space bisimilar to the
  given finite topological space. If the class of topological spaces
  in question is closed under refinement of topology, we can do it the
  same way as we have done in the Theorem~\ref{thm:t2}.
\end{proof}

\section{Hybrid logic}

In this section we will denote by $\Log(K)$ a set of formulas in the
hybrid language $H(E)$ (with nominals, $@$ and universal modality)
which are valid on all topological spaces in the class $K$.

In the subsequent subsections we will prove decidability of logics of
different separation axioms. Our main tool will be the notion of
topological filtration, which allows to present the information
relevant for the satisfiability of a formula in a finite structure.

\begin{defn}[Topological filtration]
  Let $\Sigma$ be a subformula-closed set of formulas and ${\mathfrak
    M}=(T, \tau, V)$ be a topological model. Define an equivalence
  relation $\eq_\Sigma$ on $T$ as follows:
$$
 w \eq_\Sigma v \textrm{ iff }\\
 \forall \varphi \in \Sigma\
 \mathfrak{M}, w \models \varphi \textrm{ iff } \mathfrak{M}, v
 \models \varphi
$$

A {\em filtration} of $\mathfrak M$ through $\Sigma$ is a model
${\mathfrak N}=(S,\sigma,W)$, defined as follows. Let $S=T/\eq_\Sigma$
and let us denote by $[s]$ an equivalence class of $\eq_\Sigma$ with
a representative $s$. For a formula $\varphi\in\Sigma$ define
$$
\Ext{\mathfrak N}{\varphi}=\{[x] \mid \mathfrak \mathfrak{M},x\models
\varphi\}
$$  
and $W(p)=\Ext{\mathfrak N}{p}$. This is well-defined, because points
from the same equivalence class satisfy the same formulas from
$\Sigma$. 

Let $\pi$ be a natural projection map $t \mapsto [t]$. Define $\sigma$
to be the finest topology that makes $\pi$ continuous (that is,
$\sigma$ is the quotient topology). 

Note that if $\Sigma=Cl(\varphi)$ (all subformulas of a single formula
$\varphi$), then any filtration by $\Sigma$ is finite (there is only
finite number of subsets of $Cl(\varphi)$).
\end{defn}


  

\subsection{$T_1$ spaces}

The class $T_1$ does not have a finite model property with respect to
the class of $T_1$ spaces: for example, the formula $i \to \Box i$ can
only be falsified on an infinite model with $T_1$ topology. In order
to prove decidability of $\Log(T_1)$ we will introduce a class of
finite topological models and prove that $T_1$ has the finite model
property with respect to that class. Then we will show that for any
formula the number of possible models from that class is bounded and
that will imply decidability.

In fact, decidability of $\Log(T_1)$ follows from the decidability of
the logic of $T_1$ spaces for $H(E)$ extended with downarrow operator
(\cite{tengasus}, section~5.4). We can justify ourselves by the fact
that the proof presented here provides us with concrete structures
that represent what $T_1$ spaces are ``from the point of view of
hybrid logic'' and that will help us later to prove complexity
results.

\begin{defn}[Finite representation of a $T_1$ model]
  \label{def:quasi-model-t1}
  
  A $T_1$ model is called {\em finitely representable} if it is
  topobisimilar to a finite topological model. A finite topological
  model where the complement of any point named by a nominal is open
  is called {\em finite representation of a $T_1$ model} (we will say
  simply {\em finite representation}, when there is no confusion).
\end{defn}

\begin{thm} \label{thm:t1_hyb}
  A formula $\varphi$ has a $T_1$ model iff it has a finitely
  representable $T_1$ model.
\end{thm}

\begin{proof}
  The left-to-right direction is proved using filtrations. Indeed, a
  filtration of any $T_1$ space is a representation of a $T_1$ model,
  as follows from the fact known from general topology (see
  \cite{gentop}) that the natural projection is open. Then it is left
  to apply the standard argument that filtration preserves
  satisfiability of all subformulas of $\varphi$.

  To prove right-to-left direction we will construct a $T_1$ model
  $\mathfrak{M}=(S,\sigma,W)$ such that the given finite
  representation $(T, \tau, V)$ is an interior image of
  $\mathfrak{M}$.

  We identify $T$ with the set of natural numbers
  $\{1,\ldots,n\}$. Suppose there are $m$ points $t_1,\ldots,t_m \in
  T$ named by a nominal. If $m=n$ then every point is named by a
  nominal and should be represented by a singleton. In this case the
  model construction process described below will produce a model with
  a finite discrete submodel.

  Let the support of $\mathfrak{M}$ be $\N$, the set of natural
  numbers. Denote $X_i=\{k\}$ for $i=t_k, 1 \leq k \leq m$ and let
  $X_i$ for $i \in T - {t_1,\ldots,t_m}$ form a partition of $\N -
  \{1,\ldots,m\}$ such that every $X_i$ is an infinite coinfinite
  set. Let $\sigma_0$ be a collection of cofinite subsets of $\N$ and
  for any subset $O \subseteq T$ denote

  $$
  \tilde{O}=\bigcup_{i \in O} X_i
  $$

  Note that $X_k$ for $m+1 \leq k \leq n$ are dense in $S$ in the
  topology $\sigma_0$. Then define the topology $\sigma$ on $\N$ to be
  generated by the following set:

  $$
  \sigma_0 \cup \{ \tilde{O} \mid O \in \tau \}
  $$

  Note that by construction of $\sigma$ if $\tilde{O}$ is open then
  $O$ is open.
  
  The valuation is defined as follows:

  $$
  W(p)=\widetilde{V(p)} \textrm{ for all } p\in \textsc{Prop} \cup
  \textsc{Nom}
  $$

  The definition of topology and valuation looks similar to the
  definition in Section~\ref{sect:bml} (indeed, the only real
  difference is in the definition of $X_k$); however, the proof works
  differently because of nominals.

  Define $f: S \to T$ to be the map that maps $X_k$ to $k$ for all $k
  \in T$. We will prove that $f$ is an interior map and that its graph
  is a total hybrid topobisimulation.

  Indeed, take an arbitrary open from $\sigma$, it will have the form
  $\tilde{O} \cup U$ where $O \in \tau$ and $U \in \sigma_0$. It can
  be represented as a union $\cup_{k \in O} (X_k \cap U)$. Here
  certain $X_k$-s are dense in $S$ in the topology $\sigma_0$, in this
  case $X_k \cap U$ is non-empty. All the other $X_k$-s are singletons
  that correspond to points named by nominals and in that case either
  $X_k \cap U$ is just $X_k$ or the intersection is empty. Denote by
  $F$ the set of those $k$-s that have an empty intersection with
  $U$. Then $f(\tilde{O} \cap U) = O \setminus F = O \cap (T \setminus
  F)$. The set $T \setminus F$ is open, because it is a complement of
  a set of points named by nominal, hence $O \setminus F$ is open. We
  have proved that $f$ is open.

  The continuity of $f$ follows easily from its definition and
  construction of $\sigma$: indeed, $f^{-1}(O)=\tilde{O}$, and if $O$
  is open, then $\tilde{O}$ is open too.

  It is easy to see that the graph $R$ of $f$ is a topobisimulation:
  it satisfies the {\bf Zig} and {\bf Zag} conditions, because it is a
  graph of a topobisimulation, it satisfies {\bf Prop} by construction
  of the valuation on $S$. It is also total and connects all points
  named by the same nominal in two models. In other words $R$ is a
  total hybrid topobisimulation.

  Now, if $(T,\tau,V),k \models \varphi$ then $\mathfrak{M},v \models
  \varphi$ for all $v \in X_k$, which finishes the proof of the
  right-to-left direction.
\end{proof}

Note that the size of a filtration through $\varphi$ is bounded by
$2^{|Cl(\varphi)|}$, hence we have an upper bound on the size of
finite representations of $T_1$ models necessary to refute
non-theorems of $Log(T_1)$. This allows us to deduce

\begin{thm}
$Log(T_1)$ is decidable.
\end{thm}

\subsection{$\Log(T_1)=\Log(T_2)$}

In Theorem~\ref{thm:t2} we used the construction of
Theorem~\ref{thm:t1} and replaced the naturals with cofinite topology
with an $n$-resolvable $T_2$ space whose existence is guaranteed by
the theorem of Feng and Masaveu. In a similar fashion we are going to
reuse the notion of finite representation of a $T_1$ model and modify
the construction of Theorem~\ref{thm:t1_hyb} in order to construct
$T_2$ models out of $T_1$ finite representations.

\begin{thm} \label{thm:t2_hyb}
  A formula $\varphi$ is satisfiable on a $T_2$ space iff there exists
  a finite representation of a $T_1$ model where $\varphi$ is
  satisfiable.
\end{thm}

\begin{proof}
 Since all $T_2$ spaces are $T_1$, the same filtration argument as in
 Theorem~\ref{thm:t1_hyb} applies here.

 Now suppose we are given a finite representation $\mathfrak{M}=(T,
 \tau, V)$ such that $\varphi$ is satisfiable on it. Let $(S,
 \sigma_0)$ be an $(n-m)$-resolvable $T_2$ space, where $n=|T|$ and
 $m$ is the number of points in the finite representation named by a
 nominal. Let $X_1', \ldots, X_{n-m}'$ be the dense subsets of $S$
 which form the partition of $S$. Note that if $n > 1$ then these sets
 have empty interiors, because if one of them doesn't then no other
 can be dense. Let $X_{n-m+1}, \ldots, X_n$ be arbitrary singleton
 subsets of $S$.  Finally, denote

 $$
 X_i = X_i' \setminus \bigcup_{j = n-m+1}^{n} X_j, \textrm{ for } 1
 \leq i \leq 1,n-m
 $$

 Since $S$ is a $T_1$ space, $X_1, \ldots, X_{n-m}$ are still dense in
 $S$.

 As usual, denote 
 
 $$
 \tilde{O}=\bigcup_{i \in O} X_i
 $$
 
 and consider a new topology $\sigma$ on $S$ generated by 

 $$
 \sigma_0 \cup \{ \tilde{O} \mid O \in \tau \}
 $$

 and the valuation 

 $$
 W(p)=\widetilde{V(p)} \textrm{ for all } p\in \textsc{Prop} \cup
 \textsc{Nom}
 $$

It is left to prove that this construction preserves satisfiability of
subformulas of $\varphi$. 

We use the same argument as in the proof of Theorem~\ref{thm:t1_hyb}
here. Indeed, we are in the same setting: $X_k$ form a partition of
$S$, some of them are singleton sets (named by nominals), others are
dense in $T$ in $\sigma_0$. Consider the map $f: S \to T$ that maps
$X_k$ to $k$ for all $k \in T$. It is continuous by its construction:
the preimage of an open $O$ is $\tilde{O}$ which is open. It is open,
because, like in Theorem~\ref{thm:t1_hyb} the image of any open
$\tilde{O} \cap U$ from $\sigma$ is $O \setminus F$ where $F$ is a set
of points named by nominals, and since that $T \setminus F$ is open
follows definition of a finite representation of the $T_1$ model, we
conclude that $f$ maps opens to opens. The graph of $f$ is a hybrid
total bisimulation, which means that it preserves satisfiability of
$H(E)$ formulas and the statement of the theorem follows.

\end{proof}

Since every $T_2$ space is a $T_1$ space, we get the following
corollary

\begin{thm}
  The logic of $T_2$ spaces coincides with the logic of $T_1$ spaces
  (and hence, is decidable).
\end{thm}

\subsection{$T_0$ spaces}

In this section we will use a similar technique to prove one more
representation/decidability result, this time for $T_0$ spaces.

\begin{prop}
  An Alexandroff space corresponding to a partial order by the
  Proposition~\ref{prop:transform} is $T_0$ and the frame that
  corresponds to a $T_0$ Alexandroff space is a partial order.
\end{prop}

\begin{proof}
  This is an easy consequence of Proposition~\ref{prop:transform}.

\end{proof}

By the Proposition above, every $T_0$ validity is a partial order
validity. The converse is not true.

Consider the countable topological space $(\N, \sigma)$ with cofinite
topology. Construct a topological space $(T, \tau)$ as follows: let
$T=\{*\} \cup \N$ and $\tau= \{ U = \{*\} \cup O \mid O \in \sigma
\}$. This is a $T_0$ space. Now, introduce a valuation that names $*$
with a nominal $i$ and consider a formula $\varphi=\Diamond (\neg i
\land \Diamond i)$. This formula is satisfied at $*$, but it is not
satisfiable on any partial order. Hence $\Log(T_0)$ is a strict subset
of the logic of partial order.

Although the counterexample we have just mentioned tells us that
$Log(T_0)$ is more complicated than the logic of partial orders, it
will serve us as the source of ideas on how one might build a $T_0$
model out of a quasi-model. We will need a different notion of a
finite representation of a model than one for $T_1$ and $T_2$ spaces
(otherwise $Log(T_0)$ would coincide with $Log(T_1)$ which is
impossible).

\begin{defn}[finite representation of a $T_0$ model]
  \label{def:quasi-model-t0}

  A {\em finitely representable $T_0$ model} is a $T_0$ topological
  model which is topobisimilar to a finite topological model. A finite
  topological model is called a {\em finite representation of a $T_0$
    model} if for every pair of points $x, y$ named by nominals, there
  exists an open neighborhood $O_x$ of $x$ such that $y \notin O_x$ or
  there exists an open neighborhood $O_y$ of $y$ such that $x \notin
  O_y$ (we will say simply {\em finite representation}, when there is
  no confusion).
\end{defn}

Once again we will describe a way to construct a topological space
(this time a $T_1$ space) that satisfies a given formula given a
finite representation that satisfies that formula. We will have as a
consequence a  

\begin{thm} \label{thm:t0_hyb}
$\Log(T_0)$ is decidable.
\end{thm}

\begin{proof}
  What we really prove here is that a formula has a $T_0$ model iff it
  has a finitely representable $T_0$ model.

  A filtration of a $T_0$ space through $Cl(\varphi)$ gives a finite
  representation of a $T_0$ model, because the natural projection is an
  open map. This construction preserves satisfiability by the same
  argument, as the one that was mentioned in the previous sections.

  The other direction of the proof goes as follows. Consider a finite
  representation $\mathfrak{M}=(T, \tau, V)$. We identify $T$ with
  natural numbers $1, \ldots, n$ and we will use such a numbering that
  $1, \ldots, m$ are named by a nominal. We construct a topological
  model $(S, \sigma, W)$ with a support $\{1, \ldots, m\} \cup \N$ and
  topology and valuation defined below. We will suppose further that $n
  \neq m$ since otherwise the finite representation is already a real
  $T_0$ that satisfies $\varphi$.

  Partition $S$ into sets $X_1, \ldots, X_n$: let $X_k=\{k\}$ for
  $1 \leq k \leq m$ and let $X_{m+1}, \ldots, X_n$ be the sets of the form
  $\{k + j (n-m) \mid 0 \leq j < \infty\}$ for $m+1 \leq k \leq n\}$.
  
  As usual, denote 
  
  $$
  \tilde{O}=\bigcup_{i \in O} X_i
  $$

  for $O \subseteq T$. Define the topology $\sigma$ to be 
  
  $$
  \{\tilde{O} \setminus F \mid O \in \tau, F \subseteq \N \textrm{ finite } \}
  $$

  Valuation is also defined in a usual way
  
  $$
  W(p)=\widetilde{V(p)} \textrm{ for all } p\in \textsc{Prop} \cup  \textsc{Nom}
  $$

  The model thus constructed is $T_0$. Any point $x$ from $\N$ can be
  separated from any other point by a set $S - \{x\}$. Since two points
  named by nominal can be separated by an open $O$ in the finite
  representation, $\tilde{O}$ will separate them in $\mathfrak M$ (that
  is where we use the fact that $T$ is a finite representation of a
  $T_0$ model).

  Consider the map $f: S \to T$ that maps $X_k$ to $k$ for all $k \in
  T$. This is an interior map and its graph is a total hybrid
  topobisimulation. 

  Indeed, note that $X_1, \ldots, X_n$ have non-empty intersection with all the
  sets of the form $T \setminus F$, where $F \subset \N$ is
  finite. Similarly, any open from $\sigma$ can be seen as a union
  $\cup_{k \in O} (X_k \cap (T \setminus F)$ for some $O \in \tau$ and
  some fixed finite $F \subset N$. Thus the image of this set under $f$ will
  be $O$, which is open.

  The continuity of $f$ follows from its construction and definition of
  the topology $\sigma$. The remaining conditions that make the graph
  of $f$ a hybrid bisimulation can be checked straightforwardly.

  Since total hybrid bisimulations preserve satisfiability of $H(E)$
  formulas, $(S, \sigma, W)$ satisfies the same formulas as $(T, \tau,
  V)$, which finishes the proof of the theorem.

\end{proof}



\subsection{Complexity}

Now, when we know that $Log(T_0)$ and $Log(T_1)$ are decidable, the
next natural question to ask is what the complexity is. The lower
bound follows from the result of Ladner \cite{ladner} that $S4$ has a
PSPACE-complete satisfiability problem.

\begin{prop}
  $Log(T_0)$ and $Log(T_1)$ have a PSPACE-hard satisfiability problem.
\end{prop}

\hyphenation{pa-ra-met-rized}

To establish an upper bound we will present a two player game
parametrized by a formula where one of the players has a winning
strategy iff the formula is satisfied on a finite representation of a
$T_0$ or $T_1$ model (and hence, is satisfiable on a $T_0$ or $T_1$
space). The amount of information on the board at the end of any play
will be polynomial in the length of the formula. Thus, it is possible
to build a polynomial space Turing machine that decides whether the
game has a winning strategy by just repeatedly analyzing all possible
plays.

We will present a different notion of a model, equivalent to the
notions of finite representation of a $T_1$ (or $T_0$) model.

\begin{defn}[Hintikka set]
  Let $\Sigma$ be a set of formulas closed under subformulas and
  single negations. A set $A \subseteq \Sigma$ is called a {\em
    Hintikka set} if it is maximal subset satisfying the following
  conditions:

  \begin{enumerate}
  \item $\bot \notin A$
  \item if $\neg \varphi \in \Sigma$ then $\varphi \in A$ iff $\neg
    \varphi \notin A$
  \item if $\varphi \land \psi \in \Sigma$ then $\varphi \land \psi
    \in A$ iff $\varphi \in A$ and $\psi \in A$
  \end{enumerate}
\end{defn}

\begin{defn}[Quasi-model] \label{def:quasi-model2}
  Let $\varphi$ be a formula and $Cl(\varphi)$ be its subformula
  closure. A tuple $(T, \tau, \lambda)$, where $(T, \tau)$ is a finite
  topological space and $\lambda$ is a function from $T$ to
  $Cl(\varphi)$ is called a {\em quasi-model} for $\varphi$ if the
  following holds:

  \begin{enumerate}
  \item $\lambda(t)$ is a Hintikka set for any $t \in T$
  \item at least for one $t \in T$, $\varphi \in \lambda(t)$
   \item\label{def:quasi-open} for all $\Box \psi \in Cl(\varphi)$,
     $\Box \psi \in \lambda(t)$ iff there exists an open $O \ni t$
     such that $\forall s \in O\ \psi \in \lambda(s)$
  \end{enumerate}
  
  If we impose extra condition on the quasi model, we are then talking
  about $T_1$ or $T_0$ quasi-models:

  \begin{itemize} 
  \item[] {\bf ($T_1$ condition for quasi-models)} if $i \in
    \lambda(t)$ where $i$ is a nominal, then $T - \{t\}$ is open.    
  \end{itemize}

  \begin{itemize}
  \item[] {\bf ($T_0$ condition for quasi-models)} for every pair of
    points $x, y$ named by nominals, there exists an open neighborhood
    $O_x$ of $x$ such that $y \notin O_x$ or there exists an open
    neighborhood $O_y$ of $y$ such that $x \notin O_y$.
  \end{itemize}
\end{defn}

\begin{lem}
  This definition is equivalent to the notion of a finite
  representation of a model in the following sense: a formula
  $\varphi$ is satisfied on finite representation of a $T_1$ ($T_0$)
  model $(S, \sigma, W)$ iff there exists a $T_1$ ($T_0$) quasi-model
  for $\varphi$.
\end{lem}

\begin{proof}
 To prove the left-to-right direction take a given finite topological
 space $S, \sigma, W$ and define a mapping $\lambda: S \to
 Cl(\varphi)$:
$$
\lambda(x) = \{\psi \in Cl(\varphi) \mid (S, \sigma, W), x \models
\psi \}
$$
 Then $(S, \sigma, \lambda)$ is a quasi-model for $\varphi$.

 Right-to-left direction: take $(S, \sigma, \lambda)$ and define
 valuation $W$:
$$
W(p) = \{x \in S \mid p \in \lambda(x) \}
$$
Then $(S, \sigma, W)$ is a finite representation. One can prove by
induction on formula structure and using
condition~\ref{def:quasi-open} in the
definition~\ref{def:quasi-model2} that for all formulas $\psi \in
Cl(\varphi)$, $\psi \in \lambda(x)$ iff $(S, \sigma, W), x \models
\psi$.

\end{proof}

The winning strategy in the game we are about to describe contains all
the necessary information to build a quasi-model that satisfies the
formula. During each play of the game a piece of model is
constructed. Since the quasi-models are a special kind of finite
topological spaces and by Proposition~\ref{prop:transform}, finite
topological spaces can be regarded as relational structures, we will
think about the quasi-models as finite relational structures.

We will prove the upper bound for $H(E)$ outright; it is not much
harder than for $H(@)$ and the result is more general. One remark must
be made, the quasi-model for $H(E)$ should satisfy one extra
condition:

\begin{itemize}
\item[] {\bf (universal modality condition) } if $E\varphi \in \lambda(x)$ then
  there exists a point $y$ such that $\varphi \in \lambda(y)$.
\end{itemize}

\begin{thm}
$Log_{H(E)}(T_0)$ is PSPACE-complete.  
\end{thm}

\begin{proof}
  For the purposes of this proof we will consider $\Diamond$ as a
  primitive operator and $\Box\varphi$ as on abbreviation of
  $\neg\Diamond\neg\varphi$. Every subformula of the form $@_i\varphi$
  can be equivalently replaced by $E(i \land \varphi)$ so we do not
  consider $@$ either.

  Here is the description of the game for a formula $\varphi$.  There
  are two players: \Abelard (male) and \Eloise (female). \Eloise plays
  by putting Hintikka sets on the board and defining a transitive
  reflexive relation $R$ on them; \Abelard introduces challenges that
  she must meet. She starts the game by putting a set $\{X_0, \ldots,
  X_k\}$ on the board and introducing a relation $R$ an them (it will
  be updated after each move). The sets and the relation must satisfy
  the following conditions:

\newlength{\rulename}
\settowidth{\rulename}{INIT-DIAMOND}

\newlength{\ruledesc}
\setlength{\ruledesc}{\textwidth} \addtolength{\ruledesc}{-\rulename} 
	\addtolength{\ruledesc}{-6\tabcolsep}

  \begin{tabular}{p{\rulename}p{\ruledesc}}
    {\sc (root)} & $X_0$ contains $\varphi$, $k \leq |Cl(\varphi)|$,\\
    {\sc (init-nom)} & no nominal occurs in two different Hintikka sets,\\
    {\sc (init-diamond)} & for all $\Diamond\chi \in Cl(\varphi)$, if
    $RX_lX_j$ and $\Diamond\chi \notin X_l$    then $\Diamond\chi
    \notin X_j$ and $\chi \notin X_j$,\\ 
    {\sc (init-univ)} & for all $X_l$ and for all $E\chi \in
    Cl(\varphi)$, $E\chi \in X_l$ iff $\chi \in X_j$ for some
    $j$,\\ 
    {\sc (init-cycles)} & $R$ has no cycles.\\
  \end{tabular}

  If the conditions do not hold, \Eloise looses immediately.
  \Abelard's turn consists of selecting a Hintikka set $X_l$
  and picking a formula $\Diamond\psi$ out of it. \Eloise must
  meet the challenge by putting a Hintikka set $Y$ on the board, such
  that the following conditions hold:

\newlength{\rulenameb}
\settowidth{\rulenameb}{DIAMOND}

\newlength{\ruledescb}
\setlength{\ruledescb}{\textwidth} \addtolength{\ruledescb}{-\rulenameb} 
	\addtolength{\ruledescb}{-6\tabcolsep}

  \begin{tabular}{p{\rulenameb}p{\ruledescb}}
  {\sc (diamond)} & $\psi \in Y$, $RX_lY$ and for all
    $\Diamond\chi \in Cl(\varphi)$, if $\Diamond\chi \notin X_l$
    then $\Diamond\chi \notin Y$ and $\chi \notin Y$,\\
  {\sc (univ)} & for all $X_l$ and for all $E\chi \in Cl(\varphi)$,
    $E\chi \in X_l$ iff $\chi \in X_j$ for some $j$,\\
  {\sc (nom)} & if $i \in Y$ for some nominal $i$ then $Y$ is one
    of the Hintikka sets \Eloise played during the first move. If this
    is the case, the game stops and she wins (unless the next rule is
    violated, in which case she loses),\\
  {\sc (cycles)} & $R$ does not have cycles that involve Hintikka
    sets that contain nominals.\\
  \end{tabular}

  If \Eloise cannot find a $Y$ that satisfies those conditions, then
  the game stops and \Abelard wins. Otherwise, \Abelard must choose a
  formula of the form $\Diamond\psi$ from the last played set (that
  is, $Y$) and the game continues in a similar way. If \Eloise manages
  to meet all \Abelard's challenges and if he has no more challenges
  to present, she wins. This does not guarantee that the game will
  stop at some point, so we introduce an extra rule. A list of
  formulas played by \Abelard is kept, if he plays a formula the
  second time, \Eloise must respond with the same Hintikka set as she
  did when he played the formula for the first time. If her set
  satisfies the conditions from the previous paragraph, \Eloise wins;
  otherwise, she loses. In any case, the game stops immediately.

  We will now prove that \Eloise has a winning strategy in the game
  iff a formula $\varphi$ has a quasi-model.

  {\bf (left-to-right direction)} Suppose that \Eloise has a winning
  strategy in the game. We build a quasi-model $(S, \sigma, \lambda)$
  for $\varphi$ as follows. Let $S_0$ be the Hintikka sets played at
  the first move --- $\{X_0, \ldots, X_k\}$. Define sets $\{S_i\}$ by
  induction; suppose $S_i$ is defined, then $S_{i+1}$ is a copy of the
  Hintikka sets played by \Eloise in reply to \Abelard moves when he
  picks sets form $S_i$ (with an exception: we do not copy sets from
  the initial move when \Eloise plays them further in the game). Let
  $S$ be the disjoint union of $S_i$. Set $Rxy$ iff for all formulas
  $\Diamond\psi\in Cl(\varphi)$, $\Diamond\psi \notin x$ implies
  $\Diamond\psi \notin y$ and $\psi \notin y$. Note that $R$ thus
  defined coincides with $R$ defined throughout the game. Note also
  that $R$ is reflexive, transitive and contains no cycles that
  involve Hintikka sets named by nominals. Let $\sigma$ consist of all
  upward closed sets (as in Proposition~\ref{prop:transform}) and
  put $\lambda(x)=x$. The topology thus defined satisfies the $T_0$
  condition for quasi-models (if it did not then $R$ would contain
  cycles with points named by nominals). The universal modality
  condition for quasi-models is taken care of by the rules of the
  game: namely, by conditions {\sc (init-univ)} and {\sc (univ)}.
  
  It is left to prove that condition~\ref{def:quasi-open} in the
  Definition~\ref{def:quasi-open} is satisfied. Suppose that
  $\Box\psi=\neg\Diamond\neg\psi \in Cl(\varphi)$ and $\Box\psi \in
  \lambda(t)$, then $\Diamond\neg\psi \notin \lambda(t)$. Then the
  conditions {\sc (init-diamond)} and {\sc (diamond)} guarantee that
  for all $s$ in the minimal upward closed set $O \ni t$,
  $\Diamond\neg\psi \notin \lambda(s)$ hence
  $\Box\psi=\neg\Diamond\neg\psi \in \lambda(s)$. By definition of
  $\sigma$, $O$ is open.

  Suppose now that $t \in O$, $\forall s\in O\ \psi \in \lambda(s)$
  where $O$ is open, or upward closed set. We need to prove that
  $\Diamond\neg\psi \notin \lambda(t)$. We will prove it by
  contradiction: if $\Diamond\neg\psi \in \lambda(t)$ then once
  \Abelard chooses this formula, \Eloise must respond with one of the
  Hintikka sets from $O$, but if she does that she breaks {\sc
    (diamond)} (because $\neg\psi \notin s$ for all $s \in O$) and
  loses. Hence, $\Diamond\neg\psi \notin \lambda(t)$.

  We have built a quasi-model from a winning strategy of \Eloise.
  
  {\bf (right-to-left direction)} Let us prove that \Eloise can read
  her winning strategy off a quasi-model $(S, \sigma, \lambda)$ for
  $\varphi$. Let $R$ be the relation of the corresponding relational
  structure obtained by the Proposition~\ref{prop:transform}.

  During her first move \Eloise picks a point $t$ such that $\varphi
  \in \lambda(t)$, for each nominal contained in $\varphi$ she picks a
  point named by that nominal, and for each subformula of $\varphi$ of
  the form $E\psi$ she picks a point $t$ such that $\psi \in
  \lambda(t)$. This move complies with the required conditions. 

  Next, when \Abelard chooses a point $X$ and a formula
  $\Diamond\psi$, \Eloise responds with a with a maximal (with respect
  to the relation $R$ understood as order relation) successor $Y$ of
  $X$ such that $Y$ contains $\psi$. Obviously this complies with {\sc
    (diamond)}, {\sc (univ)}, {\sc (nom)} and {\sc (cycles)} rules. It
  is always possible to find a maximal successor because quasi-models
  are finite. \Eloise needs to adopt this strategy to be able to
  successfully answer with the same Hintikka set when \Abelard will
  pick formula $\Diamond\psi$ again. For suppose \Eloise played $Y$ in
  response for \Abelard's challenge $\Diamond\psi$ from $X$ and
  suppose that later \Abelard picks the same formula $\Diamond\psi$
  from a set $Z$, which is a successor of $X$. Since $RXZ$ any
  successor of $Z$ containing $\psi$ will be a maximal successor of
  $X$ containing $\psi$. So $Y$ is among successors of $Z$ and can be
  played again to fulfill the rules of the game.

\end{proof}

\begin{thm}
$Log_{H(E)}(T_0)$ is PSPACE-complete.  
\end{thm}

\begin{proof}
  The game for $T_1$ is the game for $T_0$ with the following
  modifications.  {\sc (init-cycles)} and {\sc (cycles)} conditions
  are replaced with 

\vspace{1em}

  \begin{tabular}{ll}
    {\sc (no-incoming)} & points named by nominals have no incoming
    arcs
  \end{tabular}

\vspace{1em} 

\noindent and {\sc (nom)} is dropped (it is has no effect because
  of {\sc (no-incoming}).

  We only have to prove that this new rule really correspond to $T_1$
  quasi-models, the rest is taken care of in the proof for the $T_0$
  case.

  Indeed, in the model that we build out of the \Eloise's winning
  strategy no Hintikka set that contains a nominal has an incoming arc
  (because of the {\sc (no-incoming)} rule. Then a complement of any
  such point is a union of upward closed sets, hence open.

  The converse is also true: in the relational counterpart of any $T_1$
  quasi-model no nominal-named point has an incoming arc, because
  otherwise the complement of the point would not be open. Thus, when
  we build \Eloise's strategy based on a $T_1$ quasi-model, we will
  never break {\sc (no-incoming)} rule.
\end{proof}

\section{Logics of concrete spaces}

Up to now we have dealt with logics of separation axioms, that is
logics of classes of spaces. In this section we will show that logics
of certain concrete structures coincide with logics of classes of
spaces, like the class $T_1$ spaces or the class of $T_1$ spaces
without isolated points.

In basic modal logic the logics of all separation axioms coincide as a
consequence of McKinsey-Tarski theorem. This theorem implies the
inclusion $\Log(T_n) \subset Log(\Q)=S4$ (where $\Q$ is the
topological space of rational numbers) for $n=0, 1, 2, 3, 4, 5$ which
makes collapse the inclusion chain $S4 \subset \Log(T_0) \subset
\Log(T_1) \subset \ldots \subset \Log(T_5)$. The completeness results
we are going to prove collapse a similar inclusion chain
$\Log_{H(E)}(T_1) \subset \Log_{H(E)}(T_2) \subset \ldots \subset
\Log_{H(E)}(T_5)$.

But first let us present our representation result for $T_1$ under a
slightly different angle.

\begin{prop} \label{prop:union-rooted}
  Any finite representation of a $T_1$ model is topobisimilar to a
  (topological space corresponding to) disjoint union of finite rooted
  models, such that each root is named by a nominal and no other
  points are named by nominals.
\end{prop}

\begin{proof}
  Indeed, from the relational perspective finite representations of
  $T_1$ models are models where points named by nominals have no
  incoming arcs from other points (the accessibility relation is

  reflexive so there is always an arc from a point to itself). We can
  restrain our attention to the submodel generated by the points named
  by nominals, so we suppose further that any point can be reached
  from a point named by a nominal. Call this model
  $\mathfrak{M}=(W,R,V)$. Take a point $w_1$ named by a nominal and
  take its minimal successor $v_1$ which has a predecessor $w_2$ named
  by another nominal. Let the domain of the generated submodel of
  $v_1$ be $W_1$. Construct a model $\mathfrak{M}_1=((W\setminus W_1)
  \cup W_1' \cup W_1'', R', V')$, where $W_1'$ and $W_1''$ are copies
  of $W_1$. Denote by $w'$ and $w''$ members of $W_1'$ and $W_1''$
  which correspond to $w\in W_1$, then define $R'$ as follows:

$$
\begin{array}[l]{l}
Rw'v \textrm{ and } Rw''v \textrm{ for all } w \in W_1, v \in W
\textrm{ such that } Rwv\\
Rvw' \textrm{ and } Rvw' \textrm{ for all } w \in W_1, v \in W
\textrm{ such that } Rvw\\
\end{array}
$$
and define $V'$ as follows:

$$
w'' \in V'(p) \textrm{ and } w'\in V'(p) \textrm{ for all } w\in V(p)
$$

The points of the submodel of $\mathfrak{M}_1$ generated by $w_1$ have
no predecessors from the submodel generated by $w_2$. At the same time
$\mathfrak{M}_1$ is topobisimilar to $\mathfrak{M}$. By repeating the
described procedure repeatedly over $w_1$ and $\mathfrak{M}_1$ we will
finally come to a model $\mathfrak{M}_n$ such that successors of $w_1$
do not have predecessors named by nominals other than $w_1$. This way
we can ``peel off'' rooted submodels which have other points named by
nominals as their roots. At each step we have a model topobisimilar to
initial model. The model that we get in the end is a disjoint union of
rooted models, their roots named by nominals.
\end{proof}

In what follows we will call a rooted model that contains mare than
one point (the root) a {\em non-trivial} model.

\subsection{Completeness with respect to Cantor space}

The proof we present here is a slight modification of the construction
of \cite{aiello}. 

There are many equivalent definitions of Cantor space. For example,
one can regard Cantor space $\C$ as a set of paths in an infinite
binary tree, or equivalently as a set of infinite strings over
alphabet $\{0, 1\}$: 0s correspond to the path turning left, 1s
correspond to the path turning right. The topology is generated by the
basic open sets of the form:

$$
B_X = \{y \mid X \textrm{ is a prefix of } y\}
$$

\noindent where $X$ runs through all finite binary strings.

This definition is more formal than geometric, but in any case we will
not need it so much, because our completeness result is built on top
of the following result of \cite{aiello}.

\begin{defn}
  A {\em cluster} in a relational structure $(W, R)$ is a maximal set
  $A$ such that $R$ restricted to $A$ is an equivalence relation. A
  cluster is called {\em simple} if it contains only one point, and
  {\em proper} otherwise.
\end{defn}

\begin{thm}(Aiello, van Benthem, Bezhanishvili) \label{thm:cantor-modal} 

  For any topological model corresponding to a finite, rooted,
  non-trivial, reflexive and transitive model where every point except
  the root is contained in a proper cluster there exists a interior
  map from the Cantor space $\C$ onto this model.
\end{thm}

Denote $\noip$ the logic of $T_1$ spaces without isolated points. This
logic is different from $\Log(T_1)$ (consider a formula $@_i \Box i$
which is satisfiable only in a model where $i$ names an isolated
point).  We can slightly modify the proof of Theorem~\ref{thm:t1_hyb}
and prove that $\noip$ is complete with respect to finite
representations of $T_1$ models without isolated points, the
Proposition~\ref{prop:union-rooted} then is modified accordingly:

\begin{prop} \label{prop:union-rooted2} 

  Any finite representation of a $T_1$ model is topobisimilar to a
  (topological space corresponding to) disjoint union of non-trivial
  finite rooted models, such that each root is named by a nominal, no
  other points are named by nominals.
\end{prop}

The interior map that is constructed in the proof of
Theorem~\ref{thm:cantor-modal} has a nice property that it maps the
root to exactly one point of the Cantor space, that means that its
graph is a hybrid bisimulation (supposing that the root is named by a
nominal). By Proposition~\ref{prop:union-rooted} $\noip$ is complete
with respect to models which are disjoint unions of non-trivial finite
rooted transitive and reflexive submodels. This is almost the same
kind of models that we have in the Theorem~\ref{thm:cantor-modal}. We
only lack the condition that every point except the root belongs to a
proper cluster.

\begin{thm} \label{thm:cantor-complete}
  $\noip$ is complete with respect to Cantor space.
\end{thm}

\begin{proof}
  First observe that the models form the
  Proposition~\ref{prop:union-rooted} are topobisimilar to the models
  that have exactly the same properties but in addition every point
  (except points named by nominals) is a part of proper cluster. This
  can be done by taking all simple clusters and replacing them with,
  for example, two-point clusters, ensuring that all the incoming and
  outgoing edges from the initial cluster are copied to the new
  one. It is easy to see that this operation gives rise to a
  bisimulation.

  Now we have that $\Log(T_1)$ is complete with respect to models that
  are disjoint unions of models that are suitable for
  Theorem~\ref{thm:cantor-modal}. It is left to notice that Cantor
  space is homeomorphic to a disjoint union of any finite number of
  copies of itself to conclude the proof.  
\end{proof}

\begin{cor}
 Denote by $\dot{\C}$ a disjoint union of $\C$ with countably many
 isolated points. Then $\Log(T_1)=\Log(\dot{\C})$.
\end{cor}

\begin{proof}
  One can map the non-trivial rooted submodels of a finite
  representation $\mathfrak{M}$ that refutes a non-theorem of
  $\Log(T_1)$ to $\C$ as in Theorem~\ref{thm:cantor-complete}, and
  then map the isolated points of $\mathfrak{M}$ to isolated points of
  $\dot{C}$. One can do that because there is only finitely many
  isolated points in $\mathfrak{M}$. The map obtained is interior. 
\end{proof}

\begin{cor}
$\Log(\dot{\C})=\Log(T_n)$ for $n=1,2,3,4,5$.
\end{cor}

\subsection{Completeness with respect to rationals}

In this section we will prove that the logic of $T_1$ spaces without
isolated points is complete with respect to $\Q$, the rational
numbers.  The proof is inspired by the proof of completeness of $S4$
with respect to rationals from \cite{johan}. The proof uses the
following statement.

\begin{thm}[Cantor] \label{thm:cantor-order}
Every countable dense linear ordering without endpoints is isomorphic
to rational numbers.
\end{thm}

One of the consequences of Theorem~\ref{thm:cantor-order} is

\begin{prop} \label{prop:rationals-copy}
The topological space of rational numbers is homeomorphic to disjoint
union of any finite number of copies of itself.
\end{prop}

\begin{proof}
  The finite disjoint union of copies of rationals can be seen as a
  linear order which is several copies of $\Q$ (as on ordered set)
  juxtaposed. This linear order is dense, countable, without
  endpoints, hence order-isomorphic to rationals. But since the
  topology on both rationals and finite disjoint union of copies of
  $\Q$ is completely determined by the order, the two spaces are
  homeomorphic.
\end{proof}

We will turn each disjoint component of a finite representation of a
$T_1$ model into an infinite $n$-ary tree, $n\geq 2$ and show that it
is topobisimilar to rational numbers. We will then apply
Proposition~\ref{prop:rationals-copy} to conclude that there exists a
topobisimulation between any finite representation of a $T_1$ model
without isolated endpoints and $\Q$.

\begin{lem} \label{lem:tree} 

  Any finite, rooted, non-trivial, reflexive and transitive model with
  a root named by a nominal is bisimilar to the full infinite,
  reflexive and transitive $n$-ary tree, $n\geq 2$ and the root of the
  tree is mapped exactly to the root of the model by the bisimulation.
\end{lem}

\begin{proof}
  Consider a model $\mathfrak{M}$ with the root $w$ named by a
  nominal. First, remove the arc going from $w$ to itself and consider
  the unraveling of the model obtained. This will be an infinite
  irreflexive, antisymmetric, non-transitive tree $\mathfrak{M}_1$ of
  finite branching $n$ ($n$ is finite because the original model was
  finite).  It is easy to see that $\mathfrak{M}_1$ is bisimilar to
  $\mathfrak{M}$, moreover, $w$ is mapped exactly to the root of
  $\mathfrak{M}_1$ by the bisimulations.

  However, it might be the case that $\mathfrak{M}_1$ is not a full
  $n$-ary tree or it might be the case that $n=1$ ($n>0$ because
  $\mathfrak{M}$ is non-trivial). In the last case put $n=2$ as it is
  our desired number of successors for each node. Our strategy is to
  go inductively through the tree and repair it there where there is
  not enough successors. Suppose that we are standing at a point $x$
  that has $m$ successors, $m < n$. Take arbitrary successor of $x$,
  $y$ and consider a subtree $T_y$ with $y$ as root. Make $n-m$ copies
  of $T_y$ and link $x$ to their roots. By repeating this manipulation
  throughout the tree we will end up with a full infinite $n$-ary tree
  which is bisimilar to $\mathfrak{M}_1$. Consider its reflexive
  transitive closure. This tree satisfies the condition of the lemma.
  
\end{proof}

\begin{lem} \label{lem:rational-homeo} The full reflexive transitive
  infinite $n$-ary tree $\T_n$ is homeomorphic to $\Q$.
\end{lem}

\begin{proof}
  We will construct an interior map $f$ from $\T_n$ onto a countable
  dense suborder $X$ of $\Q$ which by Theorem~\ref{thm:cantor-order}
  is order-isomorphic (hence, homeomorphic) to $\Q$.

  First put $f(r)=0$ for the root $r$ of $\T_n$. Then define $f$
  inductively as follows. We will say that the root belongs to the
  level 0 of the tree, its immediate successors belong to the level 1
  etc., in general, if a point belongs to level $k$ then its immediate
  successors belong to level $k+1$. Now let $f$ be defined on a point
  $w$ of level $k$ which has successors $v_1, \ldots, v_n$. Define 

$$
\begin{array}{l}
f(v_1)=f(w)-\frac{1}{(n+1)^k},\\
f(v_m)=f(w)+\frac{l-1}{(n+1)^k}\textrm{ for } 2 \leq m \leq n\\
\end{array}
$$

The map $f$ thus constructed is an homeomorphism between $\T_n$ and
$X=f(\T_n)$. It follows from construction that it is bijective. The
proof that it is open and continuous is essentially the same as the
proof in \cite{johan}.

\end{proof}

\begin{thm}
$\noip$ is complete with respect to rational numbers. 
\end{thm}

\begin{proof}
  Follows from Proposition~\ref{prop:union-rooted2},
  Lemma~\ref{lem:tree}, Lemma~\ref{lem:rational-homeo} and
  Proposition~\ref{prop:rationals-copy}.
\end{proof}

\begin{cor}
 Denote by $\dot{\Q}$ a disjoint union of $\Q$ with countably many
 isolated points. Then $\Log(T_1)=\Log(\dot{\Q})$.
\end{cor}

\bibliography{decidability}

\end{document}